\documentclass[10pt]{article}
\usepackage{amsfonts}
\usepackage{amsmath,amscd}
\usepackage{amssymb}
\usepackage{amsthm}
\usepackage{newlfont}
\usepackage{enumerate}
\usepackage[final]{showkeys}
\newtheorem{thm}{Theorem}[section]
\newtheorem*{theorem 1}{Theorem 1}
\newtheorem*{theorem 2}{Theorem 2}

 \newtheorem{cor}[thm]{Corollary}
 
 \newtheorem{proposition}[thm]{Proposition}
 \theoremstyle{definition}
 
 \theoremstyle{remark}

\renewcommand{\leq}{\leqslant}

\newcommand{\M}{\mathcal{M}}

\begin{document}
\pagestyle{myheadings}
\markright{\hfill Finite split metacyclic...\hfill}
\title{
Finite split metacyclic groups and their 2-nilpotent multipliers}
\author{S. Aofi Al-Akbi and S. Hadi Jafari   \\ {\small Department of
Mathematics,  Mashhad Branch, Islamic Azad University,}
{\small Mashhad, Iran}
\vspace{.1cm}\\
}\maketitle

\begin{abstract}
There has been a great importance in understanding the nilpotent multipliers of finite groups in recent past. Let a group $G$ be presented as the quotient of a free group $F$ by a normal subgroup $R$. Given a positive integer $c$, the $c$-nilpotent multiplier of the group $G$ is the abelian group $\mathcal M^{(c)}(G)=(R\cap \gamma_{c+1}(F))/\gamma_{c+1}(R,F)$, where $\gamma_1(R,F)=R$, $\gamma_{c+1}(R,F)=[\gamma_c(R,F),F]$, and $\gamma_{c+1}(F)=\gamma_{c+1}(F,F)$. In particular, $\mathcal{M}^{(1)}(G)$ is the Schur multiplier of $G$. The crucial aspect of the research in to the $c$-nilpotent multipliers of groups includes either establishing their structures, or estimating their sizes and exponents. One reason for studying the $c$-nilpotent multiplier is its relevance to the isologism theory of P. Hall.
The study of Schur multiplier of finite metacyclic groups goes back to the paper by F. R. Beyl in 1973. In this article, we study the 2-nilpotent multiplier of finite split metacyclic groups with the help of their nonabelian tensor squares. In particular, we give a complete description of the triple tensor product, the triple exterior product, and the 2-nilpotent multiplier of such groups.
\\{\bf Keywords} Metacyclic group . Group action . Nonabelian tensor product . Nilpotent multiplier
\\{\bf Mathematics Subject Classification} 20F05 . 20C25 . 20J99
\end{abstract}
\maketitle
\section{\textbf{Definitions And Motivation}}
A group $G$ is said to be metacyclic if it possesses a cyclic normal subgroup $N$ such that $G/N$ is also cyclic. It is clear that a metacyclic group $G$ can be written $G=HN$ with $H\leq G,$ $N\unlhd G$ and both $H$ and $N$ cyclic. If $H\cap N=1,$ then $G$ is called the split metacyclic group. It was already known to H\"{o}lder that finite metacyclic groups can be presented on two generators and three defining relations.
When $G=HN$ is split, by taking generators $x$ of $H$ and $y$ of $N,$ it is well-known that $G$ has a presentation of the form
\begin{align*}
\langle x, y~|~x^{m}=y^{n}=1, ~xyx^{-1}=y^{l}\rangle,\tag{1}
\end{align*}
where the integers $m, n,$ and $l$ satisfying $l^{m}\equiv 1$ (mod $n$).
In the rest of the paper we will fix $m, n,$ and $l$ as the integers of the above presentation.
The purpose of this paper is to show that the $2$-nilpotent multiplier $\M^{(2)}(G)$ is the direct product of two cyclic groups of orders $(n, l-1, 1+l+...+l^{m-1})$ and $1\over n$$(n, 1+l+...+l^{m-1})(n, (l-1)^{2}),$ where by $(r, s)$ we mean the greatest common divisor of the integers $r$ and $s$.

Given a positive integer $c$, recall that the $c$-nilpotent multiplier of a group $G$ which is presented as the quotient of a free group $F$ by a normal subgroup $R$, is the abelian group $\M^{(c)}(G)=(R\cap \gamma_{c+1}(F))/\gamma_{c+1}(R,F),$ where $\gamma_1(R,F)=R$, $\gamma_{c+1}(R,F)=[\gamma_c(R,F), F]$ and $\gamma_{c+1}(F)=\gamma_{c+1}(F,F).$ It was shown by R. Baer \cite{baer} that these invariants are, up to group isomorphism, independent of the choice of free presentation of $G$. The group $\M(G)=\M^{(1)}(G)$ is more known as the Schur multiplier of $G$.

One motivation to study the $c$-nilpotent multipliers is its relevance to the isologism theory of P. Hall \cite{hall, hall1} in which groups can be classified into isologism classes (see \cite[p. 406]{burns ellis} for more details).  Another important reason is that, determining the structures of $c$-nilpotent multipliers is essential in studying generalized capability and covering groups, and would be generally useful in developing such a framework. Furthermore, although the above formula of $\M^{(c)}(G)$ looks simple enough; it is mostly a hard but interesting problem to compute it. In particular, according to the method proposed in \cite{burns ellis}, which is based on the notion of triple exterior product, we hope our different approach of studying 2-nilpotent multipliers can be extended to more general cases of c-nilpotent multipliers and provides some new insights.

Let's first recall the definition of the nonabelian tensor product of two groups $G$ and $H$, by supposing that they act on each other satisfying certain compatibility conditions \cite[pp. 178-179]{brown2}, and act on themselves by conjugation. The nonabelian tensor product $G \otimes H$, was introduced by Brown and Loday in \cite{brown1}, is defined as a group generated by the symbols $g\otimes h$, where $g\in G$ and $h \in H$, with defining relations
$$ gg' \otimes h=({^g}g'\otimes{^g}h)(g \otimes h)  \hspace{5mm} \text{and}  \hspace{5mm} g \otimes hh'=(g \otimes h)(^{h}g\otimes{^{h'}}h),$$
for all $g, g'\in G$ and $h, h'\in H$, where $^{g}g'= gg'g^{-1}$.  It was shown in \cite{ellis1} with homological methods, involving the original approach in \cite{brown2,brown1}, that $G\otimes H$ is a finite group, when $G$ and $H$ are finite groups.

Since the compatible actions are satisfied when a group acts on itself by conjugation, it leads to define the nonabelian tensor square $G\otimes G$. The nonabelian exterior square $G \wedge G$ is the quotient of $G \otimes G$ modulo the central subgroup $\nabla(G)$ generated by the elements $g\otimes g$ for all $g \in G$. The image of the generator $g\otimes g'$ in $G \wedge G$ is written $g \wedge g'$. The commutator map induces a homomorphism $[~,~]: G \wedge G \longrightarrow G$ mapping $g \wedge g'$ to $[g,g']=gg'g^{-1}g'^{-1}$, in which its kernel is isomorphic to $\M(G)$ \cite{mil}.

The first attempt at determining the tensor square of a finite split metacyclic group $(1)$ was made by Brown, Johnson, and Robertson \cite{brown2} where they achieved the favorable special case when $n$ is odd.

\begin{proposition} \label{1.1}(\cite[Proposition 15]{brown2})\label{Brown}
Let $G$ be a metacyclic group with presentation $(1)$ and $n$ is odd. Then
$G\otimes G$ is the direct product of four cyclic groups with generators
$x\otimes x,\ y\otimes y,\ (x\otimes y)(y\otimes x),\ x\otimes y,$
$ of\ orders$ $ m,\ (n,\ l-1),\ (n,\ l-1,\ 1+l+...+l^{m-1}),\ and \  (n,\ 1+l+...+l^{m-1}),$
$respectively.$
Furthermore, $\nabla(G)$=$\langle x\otimes x,\ y\otimes y,\ (x\otimes y)(y\otimes x)\rangle $  and\ $\M(G)$ is cyclic of order $\frac{1}{n}(n, l-1)(n,1+l+...+l^{m-1})$.
\end{proposition}

Next Johnson \cite{johnson} completed such a computations be evaluating $G\otimes G$ for even $n$.
\begin{proposition} \label{1.2}(\cite[Proposition 16]{johnson})\label{Brown}
Let $G$ be a metacyclic group with presentation $(1)$ and $n$ is even. Then $G\otimes G$ is the abelian group with generators $X=x\otimes x,$ $Y=y\otimes y,$ $Z=(x\otimes y)(y\otimes x),$ $T=x\otimes y,$ and $A=(y\otimes y)^{l-1},$ and relations
\begin{align*}
X^{m}&=Y^{n}=Z^{l-1}=Z^{n}=Z^{1+l+...+l^{m-1}}=T^{n}=A^{2}=A^{m}=1,
Y^{l-1}=A,\ T^{1+l+...+l^{m-1}}=A^{(l-1)(m-1)m/4}.
\end{align*}
Furthermore, $\M(G)$ is cyclic of order $\frac{1}{n}(n, l-1)(n,1+l+...+l^{m-1})$.
\end{proposition}
\noindent With the above assumptions, one could easily observe that
$\nabla(G)\cong C_{m}\times C_{(n, 2(l-1))}\times C_{(n, l-1, 1+l+...+l^{m-1})}.$

Here we need to give a brief recall of the notion of triple tensor (exterior) product of groups. The notation and terminology are the same as those in \cite{Da}. Since conjugation yields compatible actions, there is a diagonal action of $G$ on $G\otimes G$ such that $^g(g_1\otimes g_2)=~{^g}g_1 \otimes{^g}g_2$. Also the tensor square $G\otimes G$ acts on $G$ by conjugation in $G$ via the homomorphism $[~,~]: G\otimes G\longrightarrow [G,G]$ induced by the commutator map. Evidently these actions are compatible and we can thus construct the triple tensor product $\otimes ^3G=(G\otimes G)\otimes G$. Also by \cite[Definition 2.11]{brown1} the triple exterior product $\wedge^3G=(G\wedge G)\wedge G$ is obtained from the tensor product $(G\wedge G)\otimes G$ by imposing the additional relations $(g\wedge g')\otimes [g,g']=1$ for all $g, g' \in G$ (see \cite{Da} for more details). Recently, the second author in joint work with some others \cite{Fa, Da} established the explicit structures of the triple tensor (exterior) products of certain groups.

Like the above results, we'll see below that although the structures of $\wedge^3G$ and $\M^{(2)}(G)$ remain invariant for any integer $n$, but the evaluation of $\otimes^3G$ for even $n$ is somewhat different from the odd $n$. Because for even case the computations depend on the relation $T^{1+l+...+l^{m-1}} = A^{(l-1)(m-1)m/4}$ which appeared above. The main result of this paper is the following theorem.

\noindent\textbf{Main Theorem.}
Let $G$ be a finite split metacyclic group with presentation $(1)$.
\begin{flushleft}
$(i)$ If $n$ is odd, then
\\$\otimes^{3}G\cong C_{m}\times C_{(n, l-1)}\times C^{2}_{(m,n,l-1)}\times C_{(n,1+l+...+l^{m-1} )}\times C^{2}_{(n,l-1,1+l+...+l^{m-1})}\times  C_{(m,n,l-1,1+l+...+l^{m-1} )}$,
\\$(ii)$ If $n$ is even, then under the assumptions of Proposition \ref{1.2},
\\$ \otimes^{3}G $ $\cong (A\otimes G^{ab})\times(B\otimes G)$ in which $A=\langle X,Z\rangle$ and $B=\langle Y, T\rangle$ are subgroups of $G\otimes G$ such that $A\otimes G^{ab}\cong C_{m}\times C_{(m,n,l-1 )}\times C_{(n,l-1,1+l+...+l^{m-1})}\times C_{(m,n,l-1,1+l+...+l^{m-1})}$, and $B\otimes G$ is the abelian group generated by  $Y_{x}=Y\otimes x,$ $Y_{y}=Y\otimes y,$ $T_{x}=T\otimes x,$ and $T_{y}=T\otimes y,$ subject to the relations
\[Y_{x}^{(m,n,2(l-1))}= Y_{y}^{(n,l-1)}= T_{y}^{(n,4(l-1), 1+l+...+l^{m-1})}=1, Y_{x}^{l-1}=T_{y}^{l^{2}-1},\]
and either
\begin{align*}
&\begin{cases}
T_{x}^{(n, 1+l+...+l^{m-1})}=1, & ~~~~ if ~~~~ m \equiv 0\ \textnormal{or}\ 1\ (\textnormal{mod}\ 4) \\
T_{x}^{n}=1, T_{x}^{1+l+...+l^{m-1}}=T_{y}^{-(l-1)^{2}}, & ~~~~ if ~~~~ m \equiv 2\ \textnormal{or}\ 3\ (\textnormal{mod}\ 4)
\end{cases}
\end{align*}
when $T^{1+l+...+l^{m-1}}=1,$ or $T_{x}^{(n,1+l+...+l^{m-1})}=1$ when $T^{1+l+...+l^{m-1}}\neq 1$,
\\$(iii)$ $\wedge^{3}G \cong C_{(n,l-1,1+l+...+l^{m-1} )}\times C_{(n,1+l+...+l^{m-1})}$,
\\$(iv)$ $\M^{(2)}(G) \cong C_{(n,l-1,1+l+...+l^{m-1} )}\times C_{(n,1+l+...+l^{m-1} )(n,(l-1)^{2})/n}$.
\end{flushleft}

It should be mentioned that the HAP package (Ellis, 2008) of GAP \cite{Gap} has methods for checking the above results in the case when $G$ is nilpotent. We applied them for several of metacyclic groups by performing the relevant computations (for more details see \cite[Section 4]{Da}). The paper is organized as follows. In Section 2, we study some basic properties of the tensor products of groups including some relations between the generators of $\otimes^3G$. Section 3 deals with the proof of our main theorem.

\section{\textit{Computing The Triple Tensor Products}}
This section contains some results to be used throughout the rest of the paper. We start with some familiar rules of nonabelian tensor product of groups.
\begin{proposition} \label{1.3}(\cite[Proposition 2]{brown2})
Let $G$ and $H$ be two groups. Then for all $g\in G$ and $h\in H$
\\$(i)$ There are homomorphism of groups $\lambda: G\otimes H\longrightarrow G$ and $ \lambda^{\prime}:G\otimes H\longrightarrow H$ such that $\lambda(g\otimes h)={g}^{h}g^{-1}$ and $ \lambda^{\prime}(g\otimes h)={^{g}hh^{-1}}.$
\\$(ii)$ $\lambda(t)\otimes h=t ^{h}t^{-1}$ and $ g\otimes \lambda^{\prime}(t)={^{g}tt^{-1}}$, and thus $\lambda(t)\otimes \lambda^{\prime}(t_{1})=[t, t_{1}]$ for all $t, t_{1} \in G\otimes H.$ Hence, $G$ acts trivially on $Ker \lambda^{\prime}$  and $H$ acts trivially on $ Ker \lambda.$
\end{proposition}
\begin{proposition} \label{1.4}(\cite[Proposition 3]{brown2})
Let $G$ and $H$ be two groups. The following relations hold in $G\otimes H$ for all $g,g' \in G$ and $h, h' \in H$ (note $[g,h]$ may be interpreted as $g^hg^{-1}\in G$ or $^ghh^{-1}\in H$):
\begin{flushleft}
$(i)~^g(g^{-1}\otimes h)=(g\otimes h)^{-1}=~^h(g\otimes h^{-1})$,\\
$(ii)~{g'}\otimes {(^ghh^{-1})}={^{g'}} (g\otimes h)(g\otimes h)^{-1}$,\\
$(iii)~(g^hg^{-1})\otimes {h'}=(g\otimes h)^{h'}(g\otimes h)^{-1}$,\\
$(iv)~(g\otimes h)(g' \otimes h')(g\otimes h)^{-1}={^{[g,h]}}(g' \otimes h')$,\\
$(v)~[g\otimes h, g'\otimes {h'}]=(g^hg^{-1})\otimes (^{g'}{h'}h'^{-1})$.
\end{flushleft}
\end{proposition}
\begin{cor}\label{cor}
For all $a,b,c,d\in G$ the following relations hold in $G\otimes G$ and $\otimes^3G$ (to simplify the notation we shall
identify $b\otimes c\otimes d$ with $(b\otimes c)\otimes d$ in $\otimes^3G$):
\begin{flushleft}
$(i)~^{a}(b\otimes c)=(a\otimes [b,c])(b\otimes c)$,\\
$(ii)~^{a}(b\otimes c\otimes d)=(b\otimes c\otimes d)(d\otimes [b,c]\otimes a)$.
\end{flushleft}
\end{cor}
In next two subsections, given the metacyclic group $G$ with presentation (1), we list some useful relations for the elements of $\otimes^{3}G$ proceeding in a series of steps. We start with the case that $n$ is odd which is relatively straightforward. Note by Proposition \ref {1.4} that $[a\otimes b\otimes c, a^{\prime}\otimes b^{\prime}\otimes c^{\prime}]=[[a, b]\otimes c, [a^{\prime}, b^{\prime}, c^{\prime}]]$, for all $a,b,c,a',b',c'\in G$. Now by setting $G=G\otimes G$ and $H=G$ in Proposition \ref{1.3}, it tells us immediately that $Ker(\lambda': \otimes^{3}G \rightarrow G)$ is central. Since $Im\lambda'=\gamma_3(G)$ is cyclic, it follows that $\otimes^{3}G$ is abelian.
\subsection{The n odd case}
Throughout this section we let $G$ be an arbitrary finite split metacyclic group as in (1) and $n$ is odd.
\begin{align*}
{^y}(x\otimes y\otimes y)={^x}(x\otimes y\otimes y)=x\otimes y\otimes y.\tag{2}
\end{align*}
\begin{proof}
Taking $G=G\otimes G$ and $H=G$ in Proposition \ref{1.3} we get the homomorphism $\lambda: \otimes^{3}G \rightarrow G\otimes G$ so that $ \lambda (a\otimes b\otimes c)=(a\otimes b)^{c} (a\otimes b)^{-1}$, for all $a, b, c \in G$. In addition, it follows from the Proposition \ref {1.4}$(iii)$ and \cite[(5.1)]{brown2} that $x\otimes y\otimes y \in Ker \lambda$. Thus Proposition \ref{1.3} concludes the result.
\end{proof}
\vspace{-1cm}
\begin{align*}
(x\otimes y\otimes y)^{l-1} = 1.\tag{3}
\end{align*}
\begin{proof}
It follows from \cite[(5.3)]{brown2} that $(x\otimes y){^x}(x\otimes y)^{-1}=(x\otimes y)^{1-l}$. According to the defining actions in $\otimes^3G$, as $^{x\otimes y}xx^{-1} = [x,y,x]= [y^{l-1} , x]= y^{-(l-1)^{2}},$ then
 \begin{align*}
(x\otimes y\otimes y)^{(l-1)^{3}} &= (x\otimes y)^{1-l} \otimes y^{-(l-1)^{2}}\\
 &= (x\otimes y)\ ^{x}(x\otimes y)^{-1}\otimes {^{(x\otimes y)}x}x^{-1}\\
 & = [x\otimes y\otimes x ,x\otimes y\otimes x ]= 1.
 \end{align*}
On the other hand from (2) and \cite[(5.3)]{brown2} we have
\begin{align*}
(x\otimes y\otimes y) &={ ^{x}(x\otimes y\otimes y)}
={^{x}(x\otimes y )} \otimes {^{x}y}
=(x\otimes y )^{l}\otimes y^{l}
=(x\otimes y\otimes y)^{l^{2}},
\end{align*}
which implies that $(x\otimes y\otimes y)^{l^{2} -1} = 1$. Hence
$$\ 1=(x\otimes y\otimes y)^{(l-1)^{3}-(l^{2}-1)}=(x\otimes y\otimes y)^{l^{3}-4l^{2}+3l}=(x\otimes y\otimes y)^{4(l-1)}.$$
Now since $n$ is odd, it finishes the proof.
\end{proof}
\vspace{-1cm}
\begin{align*}
 ^{y}(x\otimes y\otimes x)=x\otimes y\otimes x.\tag{4}
\end{align*}
\begin{proof}
From \cite[(5.3)]{brown2} and (3) we have that
\begin{align*}
^{y}(x\otimes y\otimes x) &=((x\otimes y)\otimes {^{y}x)}\\
&= ((x\otimes y)\otimes y^{1-l}x)\\
&=(x\otimes y\otimes y^{1-l}){^{y^{1-l}}(x\otimes y\otimes x)}\\
&=(x\otimes y\otimes y)^{1-l}\ ^{y^{1-l}}(x\otimes y\otimes x)\\
&={^{y^{1-l}}(x\otimes y\otimes x)},
\end{align*}
which implies $x\otimes y\otimes x$ is fixed by $y^{l}.$ In\ addition, as $l$ and $n$ are relatively prime, there exist integers $s, t$ such that $sl+tn=1$. Hence
${^y}(x\otimes y\otimes x)={^{y^{sl+tn}}(x\otimes y\otimes x)}={^{y^{sl}}(x\otimes y\otimes x)}=x\otimes y\otimes x.$
\end{proof}
\vspace{-1cm}
\begin{align*}
 ^{x}(x\otimes y\otimes x) = (x\otimes y\otimes x)^{l}.\tag{5}
\end{align*}
\begin{proof}
By applying \cite[(5.3)]{brown2} we see that $^{x}(x\otimes y\otimes x)=(x\otimes y)^{l}\otimes x$. Now the result follows by induction on $l$ together with using (4).
\end{proof}
Setting $i.j=i(1+l+l^2+...+l^{j-1})$, then
\begin{align*}
\textnormal{$(x\otimes y)^i\otimes x^j=(x\otimes y\otimes x)^{i.j}$ and $(x\otimes y)^i\otimes y^k=(x\otimes y\otimes y)^{ik}$, for any integers $i,j,k$.}\tag{6}
\end{align*}
\begin{proof}
It follows by applying \cite[(5.2)]{brown2}, (2), (4), and (5).
\end{proof}
\subsection{The n even case}
In this subsection we are concerned with the case that $G$ is a finite split metacyclic group as in (1) and $n$ is even. As we'll see, this case is slightly more complicated with respect to the odd case. For instance, the order of $y\otimes y\otimes y$ would be discussed here while in odd case it is easily determined as a direct factor of the abelian group $\nabla(G)\otimes G^{ab}$ (see Section 3 for more details).
\begin{align*}
(y\otimes y\otimes y)^{l-1}=1. \tag{7}
\end{align*}
\begin{proof}
It is readily seen that ${^x}(y\otimes y\otimes y)=y\otimes y\otimes y^l=(y\otimes y\otimes y)^{l}$. In addition, since $y\otimes y\otimes y$ belongs to $Ker\lambda$ given in the proof of (2), it follows that ${^x}(y\otimes y\otimes y)=y\otimes y\otimes y$. The result now follows easily.
\end{proof}
\vspace{-1cm}
\begin{align*}
{^y}(x\otimes y\otimes y)=x\otimes y\otimes y.\tag{8}
\end{align*}
\begin{proof}
Using \cite[(9)]{johnson} and (7),
\begin{align*}
^{y}(x\otimes y\otimes y) &= {^{y}(x\otimes y)\otimes y}\\
&=(x\otimes y)(y\otimes y )^{l-1}\otimes y\\
&={^{x\otimes y}((y\otimes y )^{l-1} \otimes y )(x\otimes y\otimes y)}\\
&=(x\otimes y\otimes y)(y\otimes y\otimes y)^{l-1}\\
&=x\otimes y\otimes y.
\end{align*}
\end{proof}
\vspace{-1cm}
\begin{align*}
\textnormal{$(x\otimes y)^i\otimes y^j=(x\otimes y\otimes y)^{ij}$, for any integers $i,j$.} \tag{9}
\end{align*}
\begin{proof}
It can be proved by induction on $i$ together with using (8).
\end{proof}
\vspace{-1cm}
\begin{align*}
{^x}(x\otimes y\otimes y)= (x\otimes y\otimes y)^{l^{2}}. \tag{10}
\end{align*}
\begin{proof}
Invoking \cite[(10)]{johnson}, (8) and (9),
\begin{align*}
{^x}(x\otimes y\otimes y)&={^{x}(x\otimes y)\otimes {^{x}y}}\\
&=(x\otimes y)^{l} (y\otimes y)^{-l(l-1)^2/ 2}\otimes y^{l}\\
&={^{(x\otimes y)^{l}}}((y\otimes y)^{-l(l-1)^2/2}\otimes y^{l} )((x\otimes y)^{l}\otimes y^{l})\\
&=(y\otimes y\otimes y)^{-l^{2}(l-1)^2/2}(x\otimes y\otimes y)^{l^{2}}.
\end{align*}
The result follows by (7).
\end{proof}
\vspace{-1cm}
\begin{align*}
(x\otimes y\otimes y)^{1+l+...+l^{m-1}} = 1. \tag{11}
\end{align*}
\begin{proof}
Clearly the assertion holds when $(x\otimes y)^{1+l+...+l^{m-1}}=1.$ If $(x\otimes y)^{1+l+...+l^{m-1}}\neq 1$, then by Proposition \ref {1.2} as $A^2=1$ we see $(x\otimes y)^{1+l+...+l^{m-1}}=(y\otimes y)^{l-1}$. Hence $(x\otimes y)^{1+l+...+l^{m-1}}\otimes y=(y\otimes y)^{l-1}\otimes y$, and consequently (7) concludes the result.
\end{proof}
\vspace{-1cm}
\begin{align*}
 (x\otimes y\otimes y)^{4(l-1)}=1. \tag{12}
\end{align*}
\begin{proof}
As $(y\otimes y)^{2(l-1)}=1$ by \cite[(15)]{johnson}, then it follows that $y$ fixes $(x\otimes y)^{2}$. So the element $(x\otimes y\otimes y)^{2}$ = $(x\otimes y)^{2}\otimes y$ belongs to $Ker\lambda,$ where $\lambda $ is the homomorphism given in the proof of (2). Now Proposition \ref{1.3}$(ii)$ implies the group $G$ acts on $(x\otimes y\otimes y)^{2}$ trivially. Consequently it follows by (10) that
$(x\otimes y\otimes y)^{2}={{^x}(x\otimes y\otimes y)^{2}}=(x\otimes y\otimes y)^{2l^{2}}.$
Therefore $(x\otimes y\otimes y)^{2(l^{2}-1)}=1$.
On the other hand, from \cite[(10)]{johnson} and Proposition \ref{1.4} we have $(x\otimes y)^{x} (x\otimes y)^{-1}$ =$(x\otimes y)^{1-l} (y\otimes y)^{(l-1)^{2}/2}$. Thus $(x\otimes y)^{1-l}=(x\otimes y)^{x}(x\otimes y)^{-1}(y\otimes y)^{-(l-1)^{2}/2}$. Hence by applying (7) we have
\begin{align*}
(x\otimes y\otimes y)^{(l-1)^{3}}&=(x\otimes y)^{1-l} \otimes y^{-(l-1)^{2}}\\
&=((x\otimes y)^{x}(x\otimes y)^{-1} (y\otimes y)^{-(l-1)^{2}/2})\otimes y^{-(l-1)^{2}}\\
&=((y\otimes y)^{-(l-1)^{2}/2}\otimes y^{-(l-1)^{2}})((x\otimes y)^{x}(x\otimes y)^{-1}\otimes y^{-(l-1)^{2}})\\
&=(y\otimes y\otimes y)^{{(l-1)^{2} \over{2}}-(l-1)^{2}}((x\otimes y)^{x} (x\otimes y)^{-1} \otimes {^{x\otimes y}x}x^{-1})\\
&=[x\otimes y\otimes x,\ x\otimes y\otimes x]=1.
\end{align*}
Combining two last equations leads to the assertion.
\end{proof}
\vspace{-1cm}
\begin{align*}
 {^y}(x\otimes y\otimes x)=(x\otimes y\otimes x)(x\otimes y\otimes y)^{l-1}.\tag{13}
\end{align*}
\begin{proof}
Using Corollary \ref{cor}, \cite[(11)]{johnson} and (7),
\begin{align*}
{^y}(x\otimes y\otimes x)&=((x\otimes [x,y]\otimes y)(x\otimes y\otimes x)\\
&=[(x\otimes y)^{l-1}(y\otimes y)^{-(l-1)^{2}(l-2)/2}\otimes y](x\otimes y\otimes x)\\
&=(x\otimes y\otimes x)(x\otimes y\otimes y)^{l-1}.
\end{align*}
\end{proof}
\vspace{-1cm}
\begin{align*}
\textnormal{$(x\otimes y)^{i}\otimes x=(x\otimes y\otimes x)^{i}(x\otimes y\otimes y)^{{i\choose 2}(l-1)^{2}}$, for any integer $i$.} \tag{14}
\end{align*}
\begin{proof}
It follows by induction on $i$, \cite[(22)]{brown2}, (9), and (13).
\end{proof}
\vspace{-1cm}
\begin{align*}
(x\otimes y\otimes y)^{l^{2}-1}=(y\otimes y\otimes x)^{l-1}. \tag{15}
\end{align*}
\begin{proof}
Use (10) and Corollary \ref{cor},
$$(x\otimes y\otimes y)^{l^{2}}={^x}(x\otimes y\otimes y)=(y\otimes[x,y]\otimes x)(x\otimes y\otimes y)=(y\otimes y\otimes x)^{l-1}(x\otimes y\otimes y).$$
\end{proof}
\vspace{-1cm}
\begin{align*}
{^x}(x\otimes y\otimes x)=(x\otimes y\otimes x)^{l}(x\otimes y\otimes y)^{(l-1)^{2}}. \tag{16}
\end{align*}
\begin{proof}
First note that $G$ acts trivially on $y\otimes y\otimes x$. It follows by Corollary \ref{cor}, \cite[(11)]{johnson}, (14), and (15) that
\begin{align*}
{^x}(x\otimes y\otimes x)&=(x\otimes [x,y]\otimes x)(x\otimes y\otimes x)\\
&=[(x\otimes y)^{l-1}(y\otimes y)^{-(l-1)^{2}(l-2)/2}\otimes x](x\otimes y\otimes x)\\
&=(y\otimes y\otimes x)^{-(l-1)^{2}(l-2)/2}((x\otimes y)^{l-1}\otimes x)(x\otimes y\otimes x)\\
&=(y\otimes y\otimes x)^{-(l-1)^{2}(l-2)/2}(x\otimes y\otimes x)^{l-1}(x\otimes y\otimes y)^{(l-2)(l-1)^{3}/2}(x\otimes y\otimes x)\\
&=(x\otimes y\otimes x)^{l}(x\otimes y\otimes y)^{-(l-2)(l-1)^{2}}.
\end{align*}
So the result is obviously true by applying the last equation used in the proof of (12).
\end{proof}
For any integers $i$ and $j$ we have:
\begin{align*}
\tag{17}
(x\otimes y)^{i}\otimes x^{j}=
\left\{\begin{array}{lll}
                (x\otimes y\otimes x)^{i.j}\ (x\otimes y\otimes y)^{{i\choose 2}.j(l-1)^{2}} & ~~~~ if ~~~~ j \equiv 0\ \textnormal{or}\ 1\ (\textnormal{mod}\ 4)\\
                (x\otimes y\otimes x)^{i.j}\ (x\otimes y\otimes y)^{({i\choose 2}.j+i)(l-1)^{2}}& ~~~~ if ~~~~ j \equiv 2\ \textnormal{or}\ 3\ (\textnormal{mod}\ 4)
            \end{array}\right.
\end{align*}
where $i.j=i(1+l+l^2+...+l^{j-1})$.
\begin{proof}
We proceed by induction on both $i$ and $j$. Assuming $i=1$, it follows by induction on $j$ together with (16) that
\begin{align*}
x\otimes y\otimes x^{j}=
\left\{\begin{array}{lll}
               (x\otimes y\otimes x)^{1.j} & ~~~~ if ~~~~ j \equiv 0\ \textnormal{or}\ 1\ (\textnormal{mod}\ 4) \\
               (x\otimes y\otimes x)^{1.j}\ (x\otimes y\otimes y)^{(l-1)^2}&~~~~ if ~~~~ j \equiv 2\ \textnormal{or}\ 3\ (\textnormal{mod}\ 4)
            \end{array}\right.
\end{align*}
It suffices to prove the second part of the assertion and the first is similar. Using (13), $l-1$ times, we have
\begin{align*}
(x\otimes y)^{i+1}\otimes x^{j}&={^{x\otimes y}}((x\otimes y)^i\otimes x^{j})(x\otimes y\otimes x^{j})\\
&={^{x\otimes y}}((x\otimes y\otimes x)^{i.j}(x\otimes y\otimes y)^{({i\choose 2}.j+i)(l-1)^{2}})(x\otimes y\otimes x)^{1.j}(x\otimes y\otimes y)^{(l-1)^2}\\
&=(x\otimes y\otimes x)^{i.j}(x\otimes y\otimes y)^{i.j(l-1)^{2}}(x\otimes y\otimes y)^{({i\choose 2}.j+i)(l-1)^{2}}(x\otimes y\otimes x)^{1.j}(x\otimes y\otimes y)^{(l-1)^2}\\
&=(x\otimes y\otimes x)^{(i+1).j}(x\otimes y\otimes y)^{({{i+1}\choose 2}.j+i+1)(l-1)^{2}}.
\end{align*}
\end{proof}

\section{\textit{Proof of The Main Theorem}}
In this section we prove the main theorem of the paper. The proof relies on the previous section together with a computation method based on the crossed pairings. There is a significant difference between the even and odd case, and w'll treat them separately beginning with the odd case.
\subsection{The n odd case}
Let $G$ be a finite split metacyclic group as in (1) and $n$ is odd. We observe from Proposition \ref{1.1} that in this case $G\otimes G$ splits as
\[G\otimes G \cong\nabla(G)\times(G\wedge G)\cong C_{m}\times C_{(n, l-1)}\times C_{(n, l-1, 1+l+...+l^{m-1})}\times C_{(n, 1+l+...+l^{m-1})},\]
where $G\wedge G$ is isomorphic to the subgroup of $G\otimes G$ generated by $\langle x\otimes y\rangle$. This is not the case in general when $n$ is even. The groups $G\wedge G$ and $G$ act on each other in such a way that $G\wedge G\cong\langle x\otimes y\rangle$ is fixed under the action of $G$. So by applying \cite[Proposition 10]{brown2} one could describe $\otimes^3G$ as follows:
$$ \ \otimes^{3}G \cong (\nabla(G)\times(G\wedge G))\otimes G\cong(\nabla(G)\otimes G)\times (G\wedge G\otimes G).$$
Now since the groups $\nabla(G)$ and $G$ act on each other trivially, Proposition 2.4 in \cite{brown1} allows us to express $\nabla(G)\otimes G$ as the tensor product $\nabla(G)\otimes G^{ab}$ of abelian groups. Therefore by Proposition \ref{1.1} we conclude that
\[\nabla(G)\otimes G^{ab}\cong C_{m}\times C_{(n,l-1)}\times C^{2}_{(m,n,l-1)} \times C_{(n,l-1,1+l+...+l^{m-1})} \times C_{(m,n,l-1,1+l+...+l^{m-1})}.\]
\indent The main task then is to determine the cyclic invariants of $G\wedge G\otimes G$. For this purpose first we expand the arbitrary element $(x\otimes y)^{i}\otimes$ $ y^{j}x^{k}$ of $G\wedge G\otimes G$ by invoking (4) and (6):
\begin{align*}
(x\otimes y)^{i}\otimes y^{j}x^{k}&=((x\otimes y)^{i}\otimes y^{j})\ {^{y^{j}}}((x\otimes y )^{i}\otimes x^{k})\\
&=(x\otimes y\otimes y )^{ij}\ {^{y^{j}}}(x\otimes y\otimes x )^{i.k}\\
&=(x\otimes y\otimes y )^{ij}(x\otimes y\otimes x )^{i.k}.
\end{align*}
This establishes what the generators of $G\wedge G\otimes G$ are. To find the orders of these generators, and any possible relations among them, by considering our analysis of $\otimes^3G$ in Subsection 2.1, we use crossed pairing into suitable cyclic groups.
Recall that for groups $G, H$ and $L$ where $G$ and $H$ acting upon each other compatibly and acting upon themselves by conjugation, a function $\Phi: G \times H \longrightarrow L$ is called a crossed pairing if
\[\Phi(gg',h)=\Phi(^{g} g', ^{g} h) \Phi (g,h) \hspace{5mm} \textnormal{and} \hspace{5mm} \Phi(g, hh')= \Phi(g,h) \Phi(^{h}g, ^{h}h'),\]
for all $g, g' \in G$, $h, h'\in H$. Clearly any crossed pairing $\Phi: G \times H \longrightarrow L$ determines a unique homomorphism $\Phi^*: G \otimes H \longrightarrow L$ such that $\Phi^*(g\otimes h)=\Phi(g,h).$

Now define $\Phi: (G\wedge G )\times G \rightarrow C_{(n , 1+l+...+l^{m-1} )} \times C_{(n ,l-1 , 1+l+...+l^{m-1} )}$ by $\Phi((x\otimes y)^{i},y^{j}x^{k})=a^{i.k}\ b^{ij}$ and let $((x\otimes y)^{i}, y^{j}x^{k})=((x\otimes y)^{i^{\prime}}, y^{j^{\prime}}x^{k^{\prime}})$. Then with setting $ r=(n, 1+l+...+l^{m-1})$, it follows that $i^{\prime}=i+ur, j^{\prime}=j+vn$, and $k^{\prime}=k+wm,$ for some integers $ u, v, w$. So $a^{i'k'}b^{i'j'}=a^{i.k}\ b^{ij}$, which implies that $\Phi$ is well-defined.
We prove that $\Phi$ is a crossed pairing. The following verifies the first rule for a crossed pairing:
\begin{align*}
\Phi(^{(x\otimes y)^{i}}(x\otimes y)^{i^{\prime}},\ ^{(x\otimes y)^{i}}y^{j}x^{k})\Phi((x\otimes y)^{i},\ y^{j}x^{k})&=\Phi((x\otimes y)^{i^{\prime}},\ y^{j+i(l-1)(1-l^{k})}x^{k})\Phi((x\otimes y)^{i},\ y^{j}x^{k})\\
&=a^{i^{\prime}.k}b^{i^{\prime}j+i^{\prime}i(l-1)(1-l^{k})}a^{i.k}\ b^{ij}\\
&= a^{(i+i^{\prime}).k}\ b^{(i+i^{\prime})j}\\
&=\Phi((x\otimes y)^{i}(x\otimes y)^{i^{\prime}}\otimes y^{j}x^{k}).
\end{align*}
The other rule follows by \cite[(22), (5.2), (5.3)]{brown2}:
\begin{align*}
\Phi((x\otimes y)^{i}, y^{j}x^{k})\Phi(^{y^{j}x^{k}}(x\otimes y)^{i}, ^{y^{j}x^{k}}y^{j^{\prime}}x^{k^{\prime}})&=\Phi((x\otimes y)^{i},\ y^{j}x^{k})\Phi((x\otimes y)^{il^{k}}, y^{j^{\prime}l^{k}+j(1-l^{k^{\prime}})}x^{k^{\prime}})\\
&=a^{i.k} b^{ij}a^{il^{k}.k^{\prime}}\ b^{il^{k}(j^{\prime}l^{k}+j(1-l^{k^{\prime}}))}\\
&= a^{i(1+l+...+l^{k-1})+i(l^{k}+l^{k+1}+...+l^{k+k^{\prime}-1})}\ b^{i(j+j^{\prime}l^{k})}\\
&=a^{i.(k+k^{\prime})}\ b^{i(j+j^{\prime}l^{k})}\\
&=\Phi((x\otimes y)^{i},\ y^{j+j^{\prime}l^{k}}x^{k+k^{\prime}})\\
&=\Phi((x\otimes y)^{i},\ y^{j}x^{k}y^{j^{\prime}}x^{k^{\prime}}).
\end{align*}
Thus $\Phi$ induces a homomorphism ${\Phi^{\ast}}: G\wedge G \otimes G \rightarrow C_{(n , 1+l+...+l^{m-1} )}$ $\times C_{(n ,l-1 , 1+l+...+l^{m-1} )},$ which shows that the generators $x\otimes y\otimes x$ and $x\otimes y\otimes y$ are independent and have the orders divided by $(n, 1+l+...+l^{m-1})$ and $(n ,l-1 ,1+l+...+l^{m-1})$, respectively. On the other hand we simply note by Proposition \ref{1.1} and (6) that $(x\otimes y\otimes x)^{(n, 1+l+...+l^{m-1})}=1$. Likewise it follows from Proposition \ref{1.1}, (2) and (3) that $(x\otimes y\otimes y)^{(n, l-1, 1+l+...+l^{m-1})}=1$. Hence $|x\otimes y\otimes x|=(n, 1+l+...+l^{m-1})$ and $|x\otimes y\otimes y|=(n ,l-1 , 1+l+...+l^{m-1})$; from which we conclude that
\[G\wedge G\otimes G\cong  C_{(n,1+l+...+l^{m-1} )}\times C_{(n,l-1,1+l+...+l^{m-1})},\]
as desired.
Consequently as $x\otimes y\otimes [x,y]=x\otimes y\otimes y^{l-1}=(x\otimes y\otimes y)^{l-1}=1$ by (2) and (3), it follows by definition that $\wedge^{3}G\cong\ G\wedge G\otimes G.$

Now we are ready to complete the proof of the main theorem in odd case by computing $\M^{(2)}(G)$. Our method relies on a tensor product approach due to Burns and Ellis \cite{burns ellis}. Let $\gamma_3^\sharp(G)$ be the quotient group $\wedge^3G/\tau(G)$ where $\tau(G)$ is the normal subgroup of $\wedge^3G$ generated by the elements
$\langle a,b,c\rangle=((a\wedge b)\wedge{^b}c)((b\wedge c)\wedge{^c}a)((c\wedge a)\wedge{^a}b),$
for all $a,b,c\in G$. By the Hall--Witt commutator identity, the homomorphism $[~,~,~]:\wedge^3G\longrightarrow G$ induces a homomorphism $[~,~,~]:\gamma_3^\sharp(G)\longrightarrow G$. It is well-known \cite[Theorem 2.9]{burns ellis} that the surjection $\ker([~,~,~]:\wedge^3G\longrightarrow G)\twoheadrightarrow \M^{(2)}(G)$ given in \cite[Theorem 2.6]{burns ellis} gives rise to the
natural isomorphism
\begin{align*}
\M^{(2)}(G)\cong \ker([~,~,~]: \gamma_3^\sharp(G)\longrightarrow G).\tag{18}
\end{align*}
So in order to describe $\M^{(2)}(G),$ first, it is required to evaluate the subgroup $\tau(G)$ of $\wedge^{3}G.$ By invoking (2) and (4) together with the initial relations of $G\wedge G$ we have
\begin{align*}
  \langle x,y,y\rangle &= (x\wedge y\wedge {^{y}y)}(y\wedge y\wedge {^{y}x)}(y\wedge x\wedge {^{x}y )}\\
 &= (x\wedge y\wedge y)(y\wedge x\wedge y^{l})\\
 &= (x\wedge y\wedge y) ((x\wedge y)^{-1}\wedge y)^{l}\\
 &= (x\wedge y\wedge y)(x\wedge y\wedge y)^{-l}\\
 &= (x\wedge y\wedge y)^{1-l}=1.
\end{align*}
Analogously (2), (3) and (4) imply that
\begin{align*}
 \langle x,y,x\rangle &= (x\wedge y\wedge {^{y}x)}(y\wedge x\wedge {^{x}x)}(x\wedge x\wedge {^{x}y)}\\
 &=(x\wedge y\wedge y^{1-l}x)((x\wedge y)^{-1}\wedge x)\\
 &=(x\wedge y\wedge y^{1-l} )\ {^{y^{1-l}}(x\wedge y\wedge x)(x\wedge y\wedge x)^{-1}}\\
 &= (x\wedge y\wedge y)^{1-l}=1.
\end{align*}
Hence $\tau (G) $ is trivial and consequently $\gamma_{3}^{\sharp}(G)$ $\cong \wedge^{3}G.$ As $\gamma_{3}(G)$ = $\langle y^{(l-1)^{2}}\rangle,$ then it readily follows by (18) that
\begin{align*}
 \M^{(2)}(G) \cong C_{(n,l-1,1+l+...+l^{m-1})} \times C_{(n,1+l+...+l^{m-1} )(n,(l-1)^{2})/n}.
\end{align*}
\subsection{The n even case}
As can be expected, the $n$ even case presents more difficulties in computation of $\otimes^{3}G $; here Proposition \ref{1.2} does not give us a splitting of $G\otimes G$ into $\nabla(G)\times(G\wedge G)$, the orders of generators are not completely determined by the analysis in Subsection 2.2, and the identities found there are a bit more complex (e.g., compare the $n$ odd case with the $n$ even case in (6) and (17)).
The reason is that for the even case the subgroup $\langle x\otimes y\rangle$ is not fixed under the action of $G$ (see \cite[(9) and (10)]{johnson}).

So in order to satisfy the hypotheses of \cite[Proposition 10] {brown2} we let $G\otimes G\cong$ $A\times B$, where $A=\langle X, Z\rangle$ and $B=\langle Y,T\rangle$ are the subgroups of $G\otimes G$ constructed by using Proposition \ref {1.2}. Thus like the discussion at the beginning of Subsection 3.1 it follows that
\begin{align*}
\otimes^{3}G &\cong(A\otimes G)\times(B\otimes G)
\cong C_{m}\times C_{(m,n,l-1 )}\times C_{(n,l-1,1+l+...+l^{m-1})}\times C_{(m,n,l-1,1+l+...+l^{m-1})}\times (B\otimes G),
\end{align*}
whence it is enough to work on the cyclic invariants of the subgroup $B\otimes G.$

Assume $Y^{i}T^{j}\in B$ and $y^{s}x^{t}\in G$ for some integers $i, j, s, t,$ and put $Y_{x}=Y\otimes x, Y_{y}=Y\otimes y, T_{x}=T\otimes x, T_{y}=T\otimes y.$ From (8), (9), (13), (17), and the argument in the proof of (2) which shows that $Y_x$ is fixed under the action of $G$, we have that
\begin{align*}
Y^{i}T^{j}\otimes y^{s}x^{t}&=(Y^{i}T^{j}\otimes y^{s})\ ^{y^{s}}(Y^{i}T^{j}\otimes x^{t})\\
&=(T^j\otimes y^{s})(Y^i\otimes y^s)\ {^{y^s}}(T^j\otimes x^t)\ {^{y^s}}(Y^i\otimes x^t)\\
&=T_y^{js}\ Y_y^{is}\ {^{y^s}}(T_x^{j.t}\ T_y^{k})\ {^{y^s}}Y_x^{it}\\
&=T_y^{js}\ Y_y^{is}\ (T_x^{j.t}\ T_y^{k+s(l-1)j.t})\ Y_x^{it}\\
&=Y_x^{it}\ Y_y^{is}\ T_x^{j.t}\ T_y^{k+s(l-1)j.t+js},
\end{align*}
where $k$ is the correspondence power of $T_y$ which already occurred in (17). This tells us that $B\otimes G$ is generated by the four elements $Y_x, Y_y, T_x,$  and $T_y$.

It remains to determine the relations among these generators. As noted above, $Y_x$ is fixed under the action of $G$. From Proposition \ref {1.2}, as $Y^n=Y^{2(l-1)}=1$, it follows that $Y_{x}^n=Y_{x}^{2(l-1)}=1$. Also we have $1=Y\otimes x^m=(Y\otimes x)^m$. So $Y_{x}^{(m,n,2(l-1))}=1$. Similarly and by (7) we obtain $Y_{y}^{(n,l-1)}=1$. Since $T^n=1$, it is readily seen by (8) that $T_y^n=1$. This together with (11) and (12) conclude that $T_{y}^{(n,4(l-1),1+l+...+l^{m-1})}=1$. Furthermore, (15) gives us $Y_{x}^{l-1}=T_{y}^{l^{2}-1}$.
Now, as the relation $T^{1+l+...+l^{m-1}}=A^{(l-1)(m-1)m/4}$ in Proposition \ref {1.2} affects our discussion, we proceed with the following two cases:
\\\textbf{Case 1)} If $T^{1+l+...+l^{m-1}}=1$.

Put $(i,j)=(1,m)$ in (17). It gives either $T_{x}^{1+l+...+l^{m-1}}=1$ or $T_{x}^{1+l+...+l^{m-1}}=T_{y}^{-(l-1)^{2}}$. On the other hand, taking $(i,j)=(n,1)$ in (17) it implies either $T_{x}^{n}=T_y^{-{n\choose 2}(l-1)^{2}}$ or $T_{x}^{n}=T_y^{-({n\choose 2}+n)(l-1)^{2}}$. By the fact $T_y^n=1$ it is now trivial that $T_{x}^{n}=1$.
So we obtain the last relation equivalent to
\begin{align*}
&\begin {cases}
T_{x}^{(n,1+l+...+l^{m-1})}=1, & ~~~~ if ~~~~ m \equiv 0\ \textnormal{or}\ 1\ (\textnormal{mod}\ 4) \\
T_{x}^{n}=1, T_{x}^{1+l+...+l^{m-1}}=T_{y}^{-(l-1)^{2}}, & ~~~~ if ~~~~ m \equiv 2\ \textnormal{or}\ 3\ (\textnormal{mod}\ 4)
\end{cases}
\end{align*}
\textbf{Case 2)} If $T^{1+l+...+l^{m-1}}\neq1$.

If $l-1$ is divisible by 4, then clearly $(l-1)(m-1)m/4$ is even. When $l-1$ isn't divisible by 4 and $m \equiv 0\ \textnormal{or}\ 1\ (\textnormal{mod}\ 4)$, then again $(l-1)(m-1)m/4$ is even. Hence $T^{1+l+...+l^{m-1}}=A^{(l-1)(m-1)m/4}=1$ and we are reduced to the first case.
So it only remains the case when $l-1$ isn't divisible by 4 and $m \equiv 2\ \textnormal{or}\ 3\ (\textnormal{mod}\ 4)$. It is obvious that $(l-1)(m-1)m/4$ is odd, whence $T^{1+l+...+l^{m-1}}=Y^{l-1}$. By applying (14) we see
\begin{align*}
Y_x^{l-1}=Y^{l-1}\otimes x=T^{1+l+...+l^{m-1}}\otimes x=T_x^{1+l+...+l^{m-1}}T_y^{{1+l+...+l^{m-1}\choose 2}(l-1)^{2}}.
\end{align*}
As $(m,n,2(l-1))$ divides $l-1,$ then $Y_x^{l-1}=1$. Therefore it follows by (11) that $T_{x}^{1+l+...+l^{m-1}}=1.$ Consequently from the argument in the first case we deduce that $T_{x}^{(n,1+l+...+l^{m-1})}=1$.

Now by imposing the relation $x\otimes y\otimes [x,y]=x\otimes y\otimes y^{l-1}=(x\otimes y\otimes y)^{l-1}=1$ to $G\wedge
 G\otimes G=\langle T_x,~T_y\rangle$, it is readily obtained that
\[\wedge^{3}G\cong C_{(n,1+l+..+l^{m-1})} \times C_{(n,l-1,1+l+...+l^{m-1})}.\]
Finally, we evaluate the subgroup $\tau(G)$ of $\wedge^{3}G.$ Likewise the odd case at the end of Subsection 3.1, we first observe that $\langle x,y,y\rangle=1$.
Also by invoking (13) we get
\begin{align*}
 \langle x,y,x\rangle &= (x\wedge y\wedge {^{y}x)}(y\wedge x\wedge {^{x}x)}(x\wedge x\wedge {^{x}y)}\\
 &=(x\wedge y\wedge y^{1-l}x)((x\wedge y)^{-1}\wedge x)\\
 &=(x\wedge y\wedge y^{1-l} )\ {^{y^{1-l}}(x\wedge y\wedge x)(x\wedge y\wedge x)^{-1}}\\
 &= (x\wedge y\wedge y)^{l(1-l)}=1.
\end{align*}
As a consequence $\tau G)\cong \langle1\rangle$ and then $\gamma_{3}^{\sharp}(G)$ $\cong \wedge^{3}G.$ As before, the isomorphism (18) now implies that
\begin{align*}
 \M^{(2)}(G) \cong C_{(n,l-1,1+l+...+l^{m-1})} \times C_{(n,1+l+...+l^{m-1} )(n,(l-1)^{2})/n}.
 \end{align*}

\section*{\textbf{Statements and Declarations}}

\noindent\textbf{Competing interests:} On behalf of all authors, the corresponding author states that there is no conflict of
interest.

\small{E-mail address: \small{saeedaofi2017@gmail.com},
s.hadi\underline{ }jafari@yahoo.com}

\end{document}